\documentstyle[12pt,amsfonts]{article}
\newcommand{\ra}{\rightarrow}

\newcommand{\hs}{\hspace{1.2em}}
\newcommand{\tty}{\hspace{0.06em}}
\newcommand{\st}[2]{\rule[#1]{0mm}{#2}} 
\newcommand{\ds}{\displaystyle}
\newtheorem{theorem}{Theorem}
\newtheorem{proposition}{Proposition}
\newtheorem{lemma}{Lemma}
\newtheorem{example}{Example}
\newcommand{\lra}{\longrightarrow}

\newcommand{\mto}{\mapsto}
\newcommand{\lmto}{\longmapsto}
\newcommand{\ot}{\otimes}
\newcommand{\op}{\oplus}
\newcommand{\cd}{\cdot}
\newcommand{\Mod}{\mbox{Mod}}
\newcommand{\End}{\mbox{End}}

\newcommand{\ve}{\varepsilon}
\newcommand{\id}{\mbox{id}}

\newcommand{\BC}{{\Bbb C}}
\newcommand{\BZ}{{\Bbb Z}}
\newcommand{\BE}{{\Bbb E}}
\newcommand{\BF}{{\Bbb F}}

\newcommand{\fg}{{\frak g}}
\newcommand{\fgz}{{\frak g}_{\overline0}}
\newcommand{\fgo}{{\frak g}_{\overline1}}
\newcommand{\fgl}{{\frak gl}}
\newcommand{\fh}{{\frak h}}
\newcommand{\fsl}{{\frak sl}}
\newcommand{\fso}{{\frak so}}
\newcommand{\fsp}{{\frak sp}}
\newcommand{\fosp}{{\frak osp}}
\newcommand{\cI}{{\cal I}}
\newcommand{\cP}{{\cal P}}
\newcommand{\A}{{\cal A}}
\newcommand{\Ac}{{\cal A}^{\circ}}
\newcommand{\Ae}{{\cal A}_e}
\newcommand{\Aec}{{\cal A}_e^{\circ}}
\newcommand{\U}{{\cal U}}
\newcommand{\Uc}{{\cal U}^{\circ}}
\newcommand{\Ue}{{\cal U}_e}
\newcommand{\Uec}{{\cal U}_e^{\circ}}
\newcommand{\K}{{\cal K}}
\newcommand{\ol}{\overline}
\newcommand{\Oo}{\overline1}
\newcommand{\Oz}{\overline0}
\newcommand{\Ot}{\overline{t}}
\newcommand{\oTheta}{\overline{\Theta}}
\newcommand{\be}{\begin{eqnarray*}}
\newcommand{\ee}{\end{eqnarray*}}
\newcommand{\bea}{\begin{eqnarray}}
\newcommand{\eea}{\end{eqnarray}}

\def\1{{\rm\hbox{1\kern-.35em\hbox{1}}}}
\newcommand{\Nn}{{\Bbb N}_n} 
\newcommand{\SLq}{\mbox{SL}_q(m|\,n)}
\newcommand{\OSPq}{\mbox{OSP}_q(2|2n)} 
\newcommand{\Uq}{U_q({\frak g})} 
\newcommand{\Uqo}{U_q({\frak g}_0)} 

\newcommand{\ODe}{\overline{\Delta}}
\newcommand{\Ove}{\overline{\varepsilon}}
\newcommand{\OS}{\overline{S}}
\newcommand{\OA}{{\overline{\cal A}}}
\newcommand{\Os}{\overline{s}}

\newcommand{\oh}{\textstyle{\frac{1}{2}}}
\newcommand{\th}{\textstyle{\frac{3}{2}}}

\begin{document}

\title{\large{\bf INVARIANT INTEGRATION ON \\
CLASSICAL AND QUANTUM LIE SUPERGROUPS}}
\author{ M. Scheunert \\
\small Physikalisches Institut, Universit\"at Bonn \\
\small Nu{\ss}allee 12, 53115 Bonn, Germany
\and
R.B. Zhang \\
\small Department of Mathematics, The University of Queensland \\
\small Brisbane, QLD 4072, Australia}
\date{}
\maketitle

\begin{abstract}
Invariant integrals on Hopf superalgebras, in particular, the classical
and quantum Lie supergroups, are studied. The uniqueness (up to scalar
multiples) of a left integral is proved, and a $\BZ_2\tty$-graded version of
Maschke's theorem is discussed. A construction of left integrals is
developed for classical and quantum Lie supergroups. Applied to several
classes of examples the construction yields the left integrals in explicit
form.
\end{abstract}

\section{Introduction}

This paper studies invariant integrals on Hopf superalgebras. We shall
focus on the Hopf superalgebras of functions on classical Lie
supergroups and their quantum counter parts, developing aspects of the
general theory of integrals on them, and also establishing an explicit
construction of such integrals.

An important feature of classical and quantum Lie superalgebras is that
their finite-dimensional representations are not completely reducible.
This imposes severe restrictions on the possible integrals on the
corresponding classical and quantum Lie supergroups. We shall
extensively investigate this fact, arriving at a result which may 
be regarded as a $\BZ_2\tty$-graded version of Maschke's theorem in an
infinite dimensional setting. 

Recall that if the dual of a given finite-dimensional Hopf algebra is
semisimple, then a generalization of Maschke's theorem (see Refs.
\cite{Swe,Abe}) applies, and the invariant integral on the Hopf algebra can be
obtained by considering a Peter-Weyl type basis of the Hopf algebra. Such a
construction of integrals fails badly in the supersymmetric setting (except
for ${\rm OSP}(1|2n)$ and ${\rm OSP}_q(1|2n)$). Here we develop an explicit
construction of integrals, which can be implemented on classical Lie
supergroups and also on type I quantum supergroups. The construction can
also be adapted to produce integrals on quantum groups at roots of unity.

The study of this paper is motivated by the great importance of the Haar
measure in the theory of locally compact Lie groups. The first place we know
of where integrals in the sense of Hopf algebra theory have shown up is
Hochschild's proof of Tannaka's duality theorem for compact groups
\cite{Hoc}. Later on they played an important role in the structure theory
of finite-dimensional Hopf algebras (for example, see
\cite{LSw, MSw, Lar, LRe}).

With the appearance of quantum groups and quantum algebras, it became
obvious that integrals have to play an important role there, too. In fact,
the quantum Haar functional is a basic tool in the $C^{\ast}$-algebra approach
to quantum groups \cite{Wor}, and it can also be used to introduce topologies
on Hopf algebras which originally are defined by purely algebraic means
\cite{DKo}. Correspondingly, there are various attempts to define integration
on quantum groups, quantum spaces and their braided generalizations (see
\cite{Chr} and the references therein).

In principle, the braided case includes Hopf superalgebras as a special
example, but it seems worthwhile to investigate the super case separately.
Needless to say, there is a huge literature dealing with the integration on
supermanifolds and supergroups, but a theory of integrals on Hopf
superalgebras seems to be missing. This will be the topic of the present
work. We hope that integrals will also prove to be useful in the further
investigation of the structure and representations of quantum supergroups,
and that our results will shed some new light on the integration over
classical, i.e., undeformed Lie supergroups.

At present, we know of only one related work \cite{Phu}. In that reference,
integrals on quantum supergroups of the special linear type are
constructed by means of the $R$-matrix formalism. However, even for the
$\SLq$ quantum supergroups the techniques used and the results derived in
that paper are totally different from those to be presented here (even
though, because of the uniqueness theorem to be proved in Section 2, the
integrals on $\SLq$ constructed here and in Ref. \cite{Phu} must be
proportional).

The organization of the paper is as follows. Section 2 develops some
general theory of integrals on Hopf superalgebras and establishes results
generalizing Maschke's theorem. Section 3 studies classical Lie supergroups.
A general construction of integrals is developed, and applied to the type I
Lie supergroups, and also the type II Lie supergroups ${\rm OSP}(1|2n)$ and
${\rm OSP}(3|2)$. Section 4 extends the results to the quantum setting,
obtaining a method for constructing integrals on quantum supergroups. As
examples, the type I quantum supergroups are studied in detail. Section 5
contains a brief discussion of our results. Finally, in the Appendix we have
collected some information on the finite dual of $U(\fgl(1))$.

We close this introduction by recalling some conventions related to
$\BZ_2\tty$-graded algebraic structures. The two elements of $\BZ_2$ are
denoted by $\Oz$ and $\Oo$. Unless stated otherwise, all gradations
considered in this work will be $\BZ_2\tty$-gradations. For any superspace,
i.e. $\BZ_2\tty$-graded vector space $V = V_{\Oz} \op V_{\Oo}$\,, we define the
gradation index $[\;\,]:\, V_{\Oz} \cup V_{\Oo} \ra \BZ_2$ by
$[x]=\alpha$ if $x \in V_{\alpha}\tty$, where $\alpha \in \BZ_2\,$. 
All algebraic notions and constructions are to be
understood in the super sense, i.e., they are assumed to be consistent
with the $\BZ_2\tty$-gradations and to include the appropriate sign factors.

\section{Integrals on Hopf superalgebras}

Let $\A$ be a Hopf superalgebra with comultiplication $\Delta$, counit $\ve$,
and antipode $S$. A left integral $\int^\ell$ on $\A$ is an element of $\A^*$,
such that
\begin{equation}
  (\id_{\A}\,\ot\int^\ell)\Delta \,=\, \1_{\A}\int^\ell. \label{intleft}
\end{equation}
Equivalently, this means that
\begin{equation}
  a^{\ast} \cd \int^\ell =\, a^{\ast}(\1_{\A}) \int^\ell, \quad
                      \forall\; a^{\ast} \in A^{\ast}, \label{intexpl}
\end{equation}
where the dot denotes the multiplication in $A^{\ast}$ deduced from
$\Delta$\,. A right integral $\int^r\in\A^*$ on $\A$ is defined by
a similar requirement
\begin{equation}
  (\int^r\!\otimes\,\,\id_{\A})\Delta \,=\, \1_{\A}\int^r.
\end{equation}

Let $\A^{\rm aop,cop}$ be the Hopf superalgebra opposite to $\A$ both
regarded as an algebra and a coalgebra. Then a linear form
$\int \in \A^{\ast}$ is a left/right integral on $\A$ if and only if it
is a right/left integral on $\A^{\rm aop,cop}$. In particular, if the
antipode $S$ of $\A$ (and hence of $\A^{\rm aop,cop}$) is invertible,
then $S^{\pm 1}$ are isomorphisms of $\A$ onto $\A^{\rm aop,cop}$, and
hence $\int$ is a left integral on $\A$ if and only if $\int S^{\pm 1}$
are right integrals on $\A$\,. Thus we only need to consider left integrals
(or right integrals).

We have the following result:

\begin{theorem}
The dimension of the space of left integrals on $\A$ is not greater than 1.
In particular, any integral on $\A$ is even or odd.
\end{theorem}

{\em Proof}: The proof is carried out by reducing the problem to the
classical non-graded case. In principle, this can be viewed as an application
of Majid's bosonization \cite{Maj}, but for the present simple case the technique
has been known for quite some time. 

For notational convenience (and for reasons that will become obvious at the
end of this proof) we define a map
 \[ \tau : \BZ_2 \times \BZ_2 \lra \BC \]
by
 \[ \tau(\alpha,\beta) = (-1)^{\alpha\beta} \;,
                             \;\quad\forall\; \alpha,\beta \in \BZ_2 \,. \]
Let $\BC\,\BZ_2$ be the group Hopf algebra of $\BZ_2$\,. The canonical basis
elements will be denoted by $g_{\alpha}$\,, $\alpha \in \BZ_2$\,. In
particular, we have
  \[ g_{\alpha} g_{\beta} = g_{\alpha + \beta} \;,
                               \quad\forall\; \alpha,\beta \in \BZ_2 \;. \]
Then
  \[ \OA = \A \ot \BC\,\BZ_2  \]
is made into a usual Hopf algebra by means of the following definitions
(where $a,b \in \A$ and $\alpha,\beta \in \BZ_2$):\\[1ex]
product:
  \[ (a \ot g_{\alpha})(b \ot g_{\beta})
                    = \tau(\alpha,[b])\, a b \ot g_{\alpha + \beta} \,, \]
coproduct:
  \[ \ODe(a \ot g_{\alpha})
 = \sum_{(a)} (a_{(1)} \ot g_{[a_{(2)}] + \alpha})
                                        \ot (a_{(2)} \ot g_{\alpha}) \,, \]
(here and in the following we use Sweedler's notation for the coproduct)
                                                                    \\[1ex]
counit:
  \[ \Ove(a \ot g_{\alpha}) = \ve(a) \,, \]
antipode:
  \[ \OS(a \ot g_{\alpha})
          = \tau([a],\alpha + [a]) \, S(a) \ot g_{-\alpha - [a]}
                                                        \;. \vspace{1ex} \]

Now let $s$ be a left integral on $\A$\,, i.e., a linear form
$s \in \A^{\ast}$ such that
  \[ (\id_\A \ot s)\Delta = \1_\A \, s \]
(according to our general conventions, $\ot$ denotes the tensor product in
the {\em graded} sense). Obviously, an element $s \in \A^{\ast}$ is a left
integral if and only if its homogeneous components are. Thus we may assume
without loss of generality that $s$ is homogeneous of degree $\sigma$\,. It
is easy to see that, notwithstanding the foregoing remark, the defining
equation for a left integral is equivalent to
  \[ \sum_{(a)} a_{(1)} s(a_{(2)}) = s(a) \, \1_\A \;,
                                             \quad\forall\; a \in \A \;, \]
i.e., it takes the same form as in the non-graded case.

Define the linear form $t_{\sigma}$ on $\BC\,\BZ_2$ by
  \[ t_{\sigma}(g_{\alpha}) = \delta_{\sigma,\,\alpha} \;,
                                     \quad\forall\; \alpha \in \BZ_2 \;. \]
Then
  \[ \Os = s \ot t_{\sigma} \]
(non-graded tensor product) is a left integral on $\OA$\,. We prove this by
showing that $\Os$ satisfies the equation analogous to that given above for
$s$\,: For all $a \in \A$ and $\alpha \in \BZ_2$\,, we have
  \[ \begin{array}{rcl}
     (\id_{\OA} \ot \Os)(\ODe(a \ot g_{\alpha}))
     &\!=\!& (\id_{\OA} \ot \Os)
             \ds{\sum_{(a)} (a_{(1)} \ot g_{[a_{(2)}] + \alpha})
                                   \ot (a_{(2)} \ot g_{\alpha}) } \\[2.8ex]
     &\!=\!& \ds{\sum_{(a)} (a_{(1)} \ot g_{[a_{(2)}] + \alpha})
                         \, s(a_{(2)})\, t_{\sigma}(g_{\alpha}) } \\[3.5ex]
     &\!=\!& \ds{\sum_{(a)} (a_{(1)} \ot g_{-\sigma + \alpha})
                         \, s(a_{(2)})\, t_{\sigma}(g_{\alpha}) } \\[3.5ex]
     &\!=\!& \ds{\Big(\sum_{(a)} a_{(1)} s(a_{(2)})\Big)
                              \ot g_0 \, t_{\sigma}(g_{\alpha}) } \\[3.5ex]
     &\!=\!& s(a) t_{\sigma}(g_{\alpha}) \, \1_\A \ot g_0         \\[1.3ex]
     &\!=\!& \Os(a \ot g_{\alpha}) \, \1_\A \ot g_0 \;,
     \end{array} \]
as required. Sullivan's theorem on the uniqueness of integrals on ordinary
(non-graded) Hopf algebras (see Refs. \cite{Sul, Abe}) now implies the
corresponding result for $\A$\,.

The reader will notice that the same proof applies to arbitrary colour Hopf
algebras (and this was the other reason to introduce the map $\tau$).
                                                               \vspace{2ex}

The uniqueness result of the theorem enables us to investigate how a left
integral behaves under ``right translations''. Thus, let $\int$ be a
non-trivial left integral on $\A$\,. We know that the linear form $\int$ is
homogeneous, let $\gamma$ be its degree. We consider the linear map
  \[ g : \A \lra \A \quad,\quad g = (\int \ot\; \id)\Delta \;. \] 
Obviously, it is homogeneous of degree $\gamma$ and not equal to zero
(otherwise, $\ve \, g = \int$ would be equal to zero). Using the
coassociativity of the coproduct, it is easy to check that
\bea  (g \ot \id) \Delta \!\!&=&\!\! \Delta\, g \label{hom} \\[0.5ex]
      (\id \ot g) \Delta \!\!&=&\!\! j \, g \;, \label{int}
\eea
where
 \[ j : \A \lra \A \ot \A \quad,\quad j(a) = \1_\A \ot a \]
is the right canonical injection of $\A$ into $\A \ot \A$\,.

Now let $a^{\ast} \in \A^{\ast}$ be an arbitrary homogeneous linear form
on $\A$\,. Equation (\ref{int}) implies that $a^{\ast} \, g$ is a left
integral on $\A$ and hence proportional to $\int$. In particular,
$a^{\ast} \, g$ vanishes on the kernel of $\int$. Since this is true for all
homogeneous elements $a^{\ast} \in \A^{\ast}$, it follows that $g$ itself
vanishes on the kernel of $\int$. Consequently, there exists a unique
element $a_0 \in \A$ such that
 \[ g(a) = \big\langle \int, a \big\rangle \, a_0 \;,
                                              \quad\forall\; a \in \A \;, \]
and $a_0$ is even. Equation (\ref{hom}) now means that
 \[ \Delta(a_0) = a_0 \ot a_0 \;. \]
Since $a_0$ is non-zero (because $g$ is non-zero), we see that $a_0$ is a
group-like element of $\A$\,. Thus we have proved the following proposition.

\begin{proposition}\label{rint}
Let $\int$ be a non-trivial left integral on a Hopf superalgebra $\A$\,.
Then there exists a unique even group-like element $a_0 \in \A$ such that
 \[ (\int \ot\; {\rm id}\,) \Delta = a_0 \int \;. \]
In particular, $\int$ is also a right integral if and only if $a_0 = \1_\A$\,.
\end{proposition} \vspace{2ex}

Let $V$ be a finite-dimensional $\BZ_2\tty$-graded right $\A$-comodule, and let
\be \omega:  V \rightarrow V \ot \A \ee
be its structure map (which, according to our general conventions, is 
supposed to be even). The antipode of $\A$ enables one to introduce a right
$\A$-comodule structure on the dual space $V^*$ of $V$, with the structure
map
\be {\bar\omega}: V^{\ast} \rightarrow V^{\ast} \ot \A \ee
uniquely defined by 
 \[  \langle v^{\ast}, w \rangle \1_\A \,=\,  
(\langle \;,\; \rangle \ot M)(\id_{V^{\ast}} \ot T \ot\id_\A) 
{\bar\omega}(v^{\ast}) \ot \omega(w),
                       \quad \forall\; v^{\ast} \in V^{\ast}, \ w \in V, \]
where $T$ is the flipping map, $M$ denotes the multiplication in $\A$
and $\langle\;,\;\rangle$ is the dual space pairing. It follows that
$\End(V) = V \ot V^{\ast}$ has a natural right $\A$-comodule structure
 \[ \delta : \End(V) \ra \End(V) \ot \A \;. \]
For later use we note that a map $g \in \End(V)$ is a comodule endomorphism
of $V$ if and only if it is even and coinvariant, i.e., it satisfies
 \[ \delta(g) \,=\, g \ot \1_\A \;. \]
If $\int$ is a left integral on $\A$\,, we define the linear map  
 \[ \Phi = (\id\ot\!\int)\,\delta\,:\,\End(V) \ra \End(V) \;. \]
Consider $\Phi(m)\in\End(V)$ for any $m\in\End(V)$. Left invariance of
$\int$ immediately leads to 
 \[ \delta(\Phi(m)) \,=\, \Phi(m) \ot \1_\A \;, \]
that is, 
\begin{lemma}
${\rm Im}\,\Phi$ is contained in the subspace of coinvariant elements of
${{\rm End}}(V)$. 
\end{lemma} 

Now we consider the case when $V$ contains a sub-comodule $V_1$\,. 
Let $P\in\End(V)$ be a projection onto $V_1$\,, i.e., ${\rm Im}\,P = V_1$
and $P^2=P$. It can be easily shown that $\Phi(P)$ satisfies
 \[ \Phi(P)V \subset V_1 \quad\mbox{and}\quad
 \Phi(P)v_1=v_1 \!\int\!\1_\A \;,\quad\forall\; v_1 \in V_1\,. \]

Suppose now that $\int \1_\A \neq 0$\,. This implies that $\int$ is
even. Thus $\Phi(P)$ is even as well, and hence it is a comodule
endomorphism of $V$. It follows that ${\rm Ker}\,\Phi(P)$ is a comodule
complement of $V_1$ in $V$. Since this holds for any finite-dimensional
right $\A$-comodule $V$ and any of its sub-comodules, we conclude that all
finite-dimensional right $\A$-comodules are completely reducible. Using
the basic fact that all finitely generated comodules are
finite-dimensional, it follows by means of standard arguments (known, for
example, from the general theory of semisimple modules over rings) that all
(not necessarily finite-dimensional) right $\A$-comodules are completely
reducible.

Conversely, let $\A$ be a Hopf superalgebra such that all right
$\A$-comodules are completely reducible. In particular, $\A$ regarded as a
right $\A$-comodule with structure map $\Delta$ is completely reducible. Let
$\A_{\,0}$ be a comodule complement of $\BC\,\1_\A$ in $\A$\,. Then any
linear form $\int^r$ on $\A$ with kernel $\A_{\,0}$ is a right(!) integral
on $\A$ such that $\int^r \1_\A \neq 0$\,. Applying the foregoing to
$\A^{\rm aop,cop}$ and $\int^r$ (which is a left integral on
$\A^{\rm aop,cop}$) we conclude that all left $\A$-comodules are completely
reducible and that $\A$ also has a left integral $\int^\ell$ such that
$\int^\ell \1_\A \neq 0$\,. It should be noted that according to Larson
\cite{Lar} analogous results hold for Hopf algebras, comodules and integrals
living in an arbitrary tensor category.

Actually, much more can be said. Let $\{V(\lambda)|\lambda\in\Lambda\}$ be a
complete representative set of all finite-dimensional right $\A$-comodules,
where $\Lambda$ is some index set. Amongst these, there is a one-dimensional
comodule, $V(0)$, say, such that under the coaction, $v \mto v \ot \1_\A$\,.
We call $V(0)$ the trivial $\A$-comodule. For each $V(\lambda)$, we choose a
basis $\{v_a^{(\lambda)} | a=1,2,\ldots,\mbox{dim}V(\lambda)\}$. Then under
the coaction of $\A$\,, we have 
 \[ \omega(v_a^{(\lambda)})
                    \,=\, \sum_b v_b^{(\lambda)} \ot t^{(\lambda)}_{ba}, \]
and the $t^{(\lambda)}_{ab}$ form a Peter-Weyl type of basis for $\A$\,. If
$\int$ denotes the linear form on $\A$ defined by
 \[ \int\!\1_\A =\, 1\;, \quad \int\! t^{(\lambda)}_{ab} =0 \;, \quad 
                                              \forall \lambda \neq 0 \;, \] 
then $\int$ is both a left and right integral on $\A$ and, obviously, it is
even.

Summarizing part of our results, we have proved the following generalization
of the well-known Maschke's theorem to the case of Hopf superalgebras (see
Refs. \cite{Swe, Lar, Abe}).

\begin{proposition}
The Hopf superalgebra $\A$ admits a left integral $\int$  with
$\int\!\1_\A \neq 0$ if and only if all right $\A$-comodules are completely
reducible. 
\end{proposition}

In the present work we are mainly interested in the case where $\A$ is a
sub-Hopf-superalgebra of the finite dual $\Uc$ of a Hopf superalgebra
$\U$. The comultiplication, counit, and antipode of $\U$ will also be denoted
by $\Delta\,$, $\ve$, and $S$, respectively. In this case, if $V$ is a right
$\A$-comodule, then $V$ also has a natural left $\U$-module structure
defined by 
 \[ x \, v \,=\, (-1)^{[x][v]} (\omega(v))(x)\;,
                                     \quad \forall x\in\U, \, v\in V \,. \]
We denote by $\U\!-\!\Mod_r$ the collection of all the left $\U$-modules 
obtained from finite-dimensional right $\A$-comodules, which forms a
monoidal category. Let us now suppose in addition that $\A$ is dense in
$\U^{\ast}$, i.e., that for any non-zero element $x \in \U$, there exists
some $a \in \A$ such that $\langle a,x \rangle \neq 0$\,. Then the
above proposition is equivalent to the following statement: 
The category $\U\!-\!\Mod_r$ is semisimple if and only if $\A$ admits a left
integral which does not vanish on the identity.  

Let us close this section by the following simple remark. As above, let $\A$
be a sub-Hopf-superalgebra of $\Uc$. The even group-like elements of $\Uc$ are
exactly the characters \mbox{of $\U$,} i.e., the superalgebra homomorphisms
of $\U$ into $\BC$\,. By convention, $\A$ always contains the unit element of
$\Uc$, i.e., the counit $\ve_{{}_{\scriptstyle \U}}$ of $\U$. This is the
so-called trivial character of $\U$. Now Proposition \ref{rint} implies the
following lemma.

\begin{lemma}\label{lr}
Suppose that $\A$ does not contain any non-trivial character of $\U$. Then
every left integral on $\A$ is also a right integral.
\end{lemma}

\section{Integrals on classical supergroups}

Let $\fg=\fg_{\Oz} \op \fg_{\Oo}$ be a finite-dimensional Lie superalgebra
\cite{Kac, Sch}, where $\fg_{\Oz}$ and $\fg_{\Oo}$ are the even and odd
subspaces respectively.  We take $\U$ to be the enveloping algebra $U(\fg)$
of $\fg$\,. $U(\fg)$ contains the enveloping algebra $U(\fg_{\Oz})$ of the
Lie subalgebra $\fg_{\Oz}$ as a subalgebra. We denote $U(\fg_{\Oz})$ by
$\Ue$\,, and let
 \[ \cI:\, \Ue \lra \U \]
be the embedding, which is a Hopf superalgebra map. It is well-known that
the dual $\cI^{\ast}$ of $\cI$ induces a Hopf superalgebra map
 \[ \cP : \Uc \lra \Uec \;, \vspace{-1.5ex} \]
which is given by
 \[ \langle \cP(a), u \rangle = \langle a\,,\,\cI(u)\rangle\;,
                           \quad\forall\; a \in \Uc \,,\; u \in \Ue \;. \]

In the present work, a Lie supergroup will be defined in terms of its Hopf
superalgebra of functions, i.e., we proceed as in the usual definition of
quantum groups \cite{RTF} or quantum supergroups \cite{Man} (for a related
treatment of supergroups, see \cite{Kos} and \cite{Bos}). More precisely,
if $\fg$ is a Lie superalgebra, the superalgebra of functions on a Lie
supergroup associated to $\fg$ will be a sub-Hopf-superalgebra $\A$ of
$\Uc = U(\fg)^{\circ}$, subject to the condition that $\A$ be dense in
$U(\fg)^{\ast}$. Actually, in our discussion of integrals, this latter
property will not be used.

Thus, let $\A$ be a sub-Hopf-superalgebra of $\Uc$. We set
 \[ \cP(\A) =  \Ae \;, \] 
which is a Hopf subalgebra of $\Uec$. Then there exist the following natural
Hopf superal\-ge\-bra maps (which are injective if $\A$ is dense in $\U^{\ast}$
and, consequently, \mbox{$\Ae$ is dense in $\Ue^{\ast}$):}
\bea
 \nu:\, U(\fg) &\lra & \Ac,\nonumber \\
             x &\lmto& \nu(x)\,,\;\; \langle\nu(x)\,,\,a\rangle
                         = (-1)^{[x][a]} \langle a \,,\, x \rangle\,,
                           \quad\forall\; a \in \A \;; \label{nu}\\[1.0ex]
 \nu_e:\, \Ue  &\lra & \Aec \,,\nonumber\\
             u &\lmto& \nu_e(u) = \tilde{u}\,,\;\;
		       \langle\tilde{u} \,,\, a_0 \rangle
                         = \langle a_0\,,\, u \rangle\,,
                          \quad\forall\; a_0 \in \Ae \;;\nonumber\\[1.0ex]
 \!\!\hat{\cI} = \nu \cI:\,
           \Ue &\lra & \Ac, \nonumber \\
             u &\lmto& \hat{u}\,, \;\;
                       \langle{ \hat u} \,,\, a \rangle
                           = \langle{ \tilde u} \,,\, \cP(a) \rangle
                               = \langle \cP(a) \,,\, u \rangle
                                   = \langle a \,,\, \cI(u) \rangle \;,
                                           \quad \forall\; a\in\A\;.\quad
\eea 

Let
 \[ \int_0 : \,\Ae \lra \BC \]
be a left integral on $\Ae$ with $\int_0 \1_{\Ae} = 1$\,. The existence of
$\int_0$ depends on properties of $\fg_{\Oz}$ and $\Ae$\,. In the case when
$\fg_{\Oz}$ is semisimple or reductive as a Lie algebra, such an $\int_0$ is
known to exist and is right invariant as well. (However, see the Appendix
about the reductive case.)

\begin{lemma} \label{into} The linear form $\ds{\int_0}\cP :\, \A \lra \BC$
is left invariant with respect to $\Ue$ in the sense that
 \[ \hat{\cI}(u) \cd (\int_0 \cP \,) \,=\, \ve(u)\int_0\cP\,,
                                              \quad\forall\; u\in \Ue \;. \]
\end{lemma}

{\em Proof}: The Lemma can be confirmed by a direct calculation. For any
$ u \in \Ue$ and $a \in \A$\,, we have
\begin{eqnarray*}
       \langle\hat{u}\cdot(\int_0\cP\,)\,,\,a \rangle
 & = & \sum_{(a)}\langle\hat{u}\,,\,a_{(1)}\rangle\int_0\cP(a_{(2)})\\
 & = & \sum_{(a)}\langle\tilde{u}\,,\,\cP(a_{(1)})\rangle\int_0\cP(a_{(2)})\\
 & = & \langle\tilde{u}\ot\!\int_0\,,\,\,\Delta\cP(a)\rangle\\[0.5ex]
 & = & \langle\tilde{u}\cdot\!\int_0,\,\cP(a)\rangle = \ve(u) \int_0\cP(a) \;.
\end{eqnarray*}

Let $J = \U \fg_{\Oz}$\,. By using the Poincar\'e, Birkhoff, Witt theorem for
Lie superalgebras \cite{Sch}, one immediately sees that

\begin{lemma} 
The subspace $J$ is a left ideal of $\U$ with finite codimension.
\end{lemma}

\noindent
Consequently, the quotient space $\U/J$ is a left $\U$-module in the standard
fashion:
 \[ x(y + J) = x y + J \;,\quad \forall\; x\in\U,\: y + J \in \U/J \,. \]
Note that this module is isomorphic to the $\U$-module induced from the
trivial $\Ue$-module. According to the usual definition, an element
$z + J \in \U/J$, with $z \in \U$\,, is said to be invariant (under the action
of $\U$) if
 \[ x(z + J) = \ve(x)z + J\,,\quad\forall\; x \in \U \;. \]
Let $z + J$ be any invariant of this type, and let $\nu(z)$ be the image of
$z$ in $\Ac$ under the natural Hopf superalgebra map (\ref{nu}). Then

\begin{theorem}\label{main}
The linear form $\int = \nu(z)\cdot\int_0\cP$ is a left integral on $\A$ and
does not depend on the choice of the representative for $z + J$. If
$z \notin J$ and if the matrix elements of the $\U$-module $\U/J$ belong
to $\A$\,, the integral $\int$ is not equal to zero.
\end{theorem}

{\em Proof}: It follows from Lemma \ref{into} that for any
$X_0 \in \fg_{\Oz}$\,, $\nu(X_0)\cdot\int_0\cP = 0$\,. As $\nu$ is an algebra
homomorphism, $\nu(y)\cdot\int_0\cP = 0$ for all $y \in J$. This proves the
second part of the theorem. Now the invariance property of $z + J$ leads to
 \[ \nu(x) \cd\! \int = \ve(x)\int, \quad \forall\; x \in \U \;. \]
This implies equation (\ref{intleft}). (Indeed, since $\A$ is contained in
$\U^{\ast}$, it is sufficient to check equation (\ref{intexpl}) for all
$a^{\ast} = \nu(x)$\,, $x \in \U$.)

To prove the last part of the theorem, we choose a homogeneous basis
$(v_i)_{1 \leq i \leq r}$ of $\U/J$ such that $v_1 = \1_\U + J$. Let $\pi$
be the representation of $\U$ in $\U/J$, and let $\pi_{i,j}$ be the matrix
elements of $\pi$ with respect to the basis $(v_i)$, i.e.,
 \[ \pi(x) v_j = \sum_{i=1}^r \pi_{i,j}(x) v_i
                   \quad\mbox{if $x \in \U$}\,, \;\, 1 \leq j \leq r \,. \]
Since $v_1$ is $\Ue$-invariant, we have
 \[ \pi_{i,1}(x) = \ve_{{}_{\scriptstyle\Ue}}(x) \delta_{i,1}
                   \quad\mbox{if $x \in \Ue$}\:, \;\, 1 \leq i \leq r \,. \]
This implies that
 \[ \int\!\pi_{i,1} \,=\, (-1)^{[v_i]} \pi_{i,1}(z)
                                               \;, \quad 1 \leq i \leq r \]
(recall that we are assuming that $\int_0 \1_{\Ae} = 1$). Since $z \notin J$
and since
 \[ z + J = \pi(z) v_1 = \sum_{i=1}^r \pi_{i,1}(z) v_i \;, \]
at least one of the matrix elements $\pi_{i,1}(z)$ must be different from
zero. This proves the theorem.

We notice that
  \[ \int \! \1_\A \,=\, \ve(z) \int_0 \1_{\Ae} \,. \]
Taking for granted that $\int_0 \1_{\Ae}$ is different from zero, we see that
$\int \1_\A \neq 0$ if and only if $\ve(z) \neq 0$\,. \vspace{1.5ex}

\noindent
{\it Remark\,\,} Suppose that the Lie algebra $\fgz$ is reductive, and that
the adjoint representation of the center of $\fgz$ in $\fgo$ is
diagonalizable. Then the subspace of $\U$-invariant elements of $\U/J$ is at
most one-dimensional. This follows at once from Theorem 1 and Theorem 2,
applied to a suitable sub-Hopf-superalgebra $\A$ of $\Uc$ (see the Appendix).
\vspace{2.0ex}

Let us now consider examples.

\begin{example} The Berezin integral\end{example} 
Consider the purely odd Lie superalgebra $\fg = \fg_{\Oo}$ with the basis
$\{\xi_i\,,\,\,i=1,2,\ldots,n\}$ and with the super bracket
 \[ [\xi_i \,,\, \xi_j] = 0 \;,\quad\forall\; i,j \;. \]
Obviously, $\U = U(\fg)$ is the Grassmann algebra on the $n$ generators 
$\xi_i$\,, and $\Ue = \BC\,\1_\U$\,. It is well-known that $\U$ has the basis
 \[ \Xi_{j_1\cdots\,j_{\ell}} = \xi_{j_1}\cdots\,\xi_{j_{\ell}} \;,
                        \quad 1 \leq j_1 < \cdots < j_{\ell} \leq n  \;, \]
where the $\ell = 0$ element is understood to be the unity. The Hopf
structure of $\U$ is the standard one for enveloping algebras of Lie
superalgebras.

Introduce a basis $\{\Theta_{i_1\cdots\,i_k} \,,\,
                 1 \leq i_1 < \cdots < i_k \leq n \,,\,\,0 \leq k \leq n\}$
for $\U^{\ast}$ (the $k = 0$ case corresponds to the unit element) such that
 \[\left\langle\Theta_{i_1\cdots\,i_k},\,\Xi_{j_1\cdots\,j_{\ell}}\right\rangle
 = (-1)^{\frac{1}{2} k(k-1)}\delta_{kl}\,\delta_{i_1j_1}\cdots\delta_{i_kj_k}
                                                                       \;, \]
and set 
 \[ \theta_i = \Theta_i \,,\quad i = 1,\,2,\,\ldots,\,n \;. \]
The cocommutativity of $\U$ implies that
 \[ \theta_i\theta_j + \theta_j\theta_i = 0 \,,\quad\forall\; i,j. \]
It is also easy to show that
 \[ \Theta_{i_1\cdots\,i_k} = \theta_{i_1}\cdots\theta_{i_k}\,,\quad
                                                    i_1 < \cdots <i_k \;. \]
As a Hopf superalgebra, $\U^{\ast}$ has the unique comultiplication such
that
 \[ \Delta(\theta_i) = \theta_i \ot \1 + \1 \ot \theta_i \;, \]
the counit is fixed by
 \[ \ve(\theta_i) = 0 \;, \]
and the antipode is specified by
 \[ S(\theta_i) = -\theta_i \;. \]

Since $(\U^{\ast})^{\ast} \cong \U$ in this case, we make the identification.
It is obvious that $\int_0\cP = \1_\U$\,. Moreover, the $\U$-invariant
elements of $\U$ are the scalar multiples of $\Xi_{1\,2\,\cdots\,n}$\,.
Thus upon choosing an appropriate normalization we obtain the unique integral
 \[ \int = (-1)^{\frac{1}{2} n(n-1)}\xi_1\xi_2\cdots\xi_n \;, \]
which yields the standard Berezin integral on the Grassmanian algebra
$\U^{\ast}$\,:
 \[ \begin{array}{rcl}
    \int\theta_{i_1}\theta_{i_2}\cdots\theta_{i_k} \!\!&=&\!\! 0\,,
                                         \quad\mbox{if $k < n$} \;, \\[1.0ex]
    \int\theta_1\theta_2\cdots\theta_n \!\!&=&\!\! 1 \;.
    \end{array} \]

To explain the left (and right) invariance of $\int$ in more familiar terms,
note that if $P(\theta)$ is any polynomial in the $\theta_i$'s, then
 \[ \Delta P(\theta) = P(\theta \ot \1 + \1 \ot \theta) \;. \]
Left invariance of the integral means
 \[ (\id \ot\! \int)\Delta(P(\theta)) = \int\! P(\theta) \;. \]
One may write $\1 \ot \theta_i$ as $\theta_i$\,, and denote $\theta_i \ot \1$
by $\lambda_i$\,, which is regarded as an independent Grassmann member. The
above equation states that
 \[ \int_\theta P(\theta + \lambda) = \int\! P(\theta)\;, \]
where the subscript $\theta$ on the left hand side indicates the fact that
the ``integration'' is carried out over the $\theta$'s. The last equation is
nothing but the translational invariance of the Berezin integral.

\begin{example} The Lie supergroup $\mbox{\rm SL}(m|\,n)$\end{example}
Let $\fg$ denote the Lie superalgebra $\fsl(m|\,n)$, which we shall regard
as a subalgebra of the general linear Lie superalgebra ${\fgl}(m|\,n)$.  
Let $\{ E_{ab} \,|\, a,b=1,2,\ldots,m+n\}$ be the standard homogeneous basis of
${\fgl}(m|\,n)$, which satisfies the commmutation relations
 \[ [E_{ab}\,,\, E_{cd}] \,=\, \delta_{bc} E_{ad}
                     - (-1)^{[E_{ab}][E_{cd}]} \delta_{da} E_{cb} \;, \]
where $[\cd\tty\,,\, \cd]$ should be understood as the graded brackets, namely,
it is symmetric when both arguments are odd, and antisymmetric otherwise.

The standard basis for $\fg$ then is given by 
 \[ E_{ab}\,,\;\; a\neq b \;; \;\; 
         h_a = E_{aa}-(-1)^{\delta_{am}}E_{a+1,\,a+1}\;, \;\; a < m+n \;.\] 
The maximal even subalgebra of $\fg$ is
$\fg_{\Oz} = {\fsl}(m) \op {\fsl}(n) \op {\fgl}(1)$. Let $\fg_{\Oo+}$ be the
odd subalgebra spanned by
$E_{i\mu}$\,,\, $i \leq m$\,, $\:\mu > m$\,, and $\fg_{\Oo-}$ be that spanned
by $E_{\mu i}$\,. Then $\fg$ is the direct sum
$\fg=\fg_{\Oo-} \op \fg_{\Oz} \op \fg_{\Oo+}$ (as vector spaces). Under the
Lie superbracket,
\bea
 [\fg_{\Oo+}\,,\,\fg_{\Oo+}] &\!=\!& \{0\} \,,\nonumber \\
 \left[\fg_{\Oo-}\,,\,\fg_{\Oo-} \right] &\!=\!& \{0\} \,,\nonumber \\
 \left[\fg_{\Oz}\,,\,\fg_{\Oo\pm}\right] &\!\subset\!& \fg_{\Oo\pm}
							     \,,\nonumber\\
 \left[\fg_{\Oo+}\,,\,\fg_{\Oo-}\right] &\!\subset &\!\fg_{\Oz} \,.
							       \label{same}
\eea

Next, we observe that $U(\fg_{\Oo+})$ and $U(\fg_{\Oo-})$ are both isomorphic
to the Grassmann algebra on $mn$ generators. The subspaces of the highest
Grassmann degree in $U(\fg_{\Oo+})$ and $U(\fg_{\Oo-})$ are both
1-dimensional. We choose the following bases for them, respectively, 
\be 
       \BE &\!=\!& \BE_{\,m}\,\BE_{\, m-1}\cdots\,\BE_{\,1}\,,\\
       \BF &\!=\!& \BF_1\,\BF_2\,\cdots\,\BF_m\,, \vspace{-2.0ex}
\ee
where \vspace{-2.0ex}
\be
 \BE_{\,i} &\!=\!& E_{i,\,m+1}\,E_{i,\,m+2}\,\cdots\, E_{i,\,m+n}\,,\\
     \BF_i &\!=\!& E_{m+n,\tty i}\,E_{m+n-1,\tty i}\,\cdots\, E_{m+1,\tty i}\,.
\ee 
Then we have \newline
\begin{minipage}[t]{6.25in} \vspace{-4ex}
 \[ [X,\,\BE\,] = [X,\,\BF\,] = 0 \,,\;\;\forall\; X \in \fg_{\Oz} 
                                                       \vspace{-1.2ex}\,,\]
\[\begin{array}{rcl}
 \xi_{+} \BE &\!=\!& 0\,,\;\;\forall\; \xi_{+} \in \fg_{\Oo+}\,,\\[0.6ex]
 \xi_{-} \BF &\!=\!& 0\,,\;\;\forall\; \xi_{-} \in \fg_{\Oo-}\,,\\[-0.3ex]
\end{array} \]
 \[ \xi_{-} \BE - (-1)^{mn} \BE\, \xi_{-} \in U(\fg)\fg_{\Oz}\,,\;\;
                          \forall\; \xi_{-} \in \fg_{\Oo-} \,.
                                                          \vspace{1.0ex} \]
\end{minipage}
Defining
 \[ \Gamma = \BE\,\BF \,, \]
it follows that
 \[ X\Gamma \in U(\fg)\fg_{\Oz}\,,\;\; \forall\; X \in \fg \,. \]

Let $t$ be the defining representation of $\fsl(m|\,n)$, with
\begin{eqnarray*}
 t(E_{ab}) &\!=\!& e_{ab}\,,\quad a\neq b \;, \\
 t(h_a) &\!=\!& e_{aa} - (-1)^{\delta_{am}} e_{a+1,\,a+1} \,,
\end{eqnarray*}
where the $e_{ab}$'s are the matrix units, and let $t_{ab}$\,,\,
$a,b=1,\,2,\,\ldots,\,m+n$\,, be the elements of $\Uc = U(\fg)^{\circ}$
defined by
 \[ \left(t_{ab}(x)\right)_{a,\,b=1}^{m+n} = t(x)\,,\;\;
                                                 \forall\; x\in U(\fg). \]
Moreover, let $\Ot$ be the dual representation of $t$, and let us similarly
introduce the matrix elements $\Ot_{ab} \in U(\fg)^{\circ}$ of $\Ot$. We
note that
 \[ \sum_c \Ot_{ca} t_{cb} (-1)^{\left([a]+[c]\right)\left([b]+\Oo\right)}
                                                    \,=\, \delta_{ab} \,. \]

The standard comultiplication on $U(\fg)$ is super cocommutative. Therefore
the finite dual $U(\fg)^{\circ}$ is a super commutative Hopf superalgebra.
The matrix elements $t_{ab}$ and ${\Ot}_{ab}$ of the vector representation
and the dual vector representation generate a sub-Hopf-superalgebra $\A$ of
$U(\fg)^{\circ}$, with the comultiplication
\be
 \Delta(t_{ab}) &\!=\!& \sum_c t_{ac}\ot t_{cb}\, 
                                           (-1)^{([a]+[ c])([c]+[b])}\;, \\
 \Delta({\Ot}_{ab}) &\!=\!& \sum_c {\Ot}_{ac}\ot {\Ot}_{cb} 
                                           (-1)^{([a]+[c])([c]+[b])} \;,  
\ee
the counit $\ve(t_{ab}) = \ve({\Ot}_{ab}) = \delta_{ab}$\,, and the
involutary antipode \vspace{0.3ex}
$S(t_{ab}) = (-1)^{[a]([a]+[b])}{\Ot}_{ba}$\,, where
$[a] = \left\{\begin{array}{ll}
          \Oz\,,&\! a \leq m\,,\\
          \Oo\,,&\! a > m\,. \end{array}\right.$  
An important fact is that  

\begin{proposition}\label{Harish}
The subspace $\A$ is dense in $U(\fsl(m|\,n))^{\ast}$.
\end{proposition}

{\em Proof}: This follows from a slight strengthening of a theorem which in
the non-graded case is due to Harish-Chandra. Let $V$ be a
finite-dimensional graded vector space and let $\fg$ be a graded subalgebra
of the Lie superalgebra $\fsl(V)$. We regard $V$ as a $\fg$-module. Arguing
as in the non-graded case (see the proof of Theorem 2.5.7 in Ref. \cite{Dix})
one can easily prove that for any non-zero element $x \in U(\fg)$ there
exists an integer $r \geq 0$ such that $x$ acts non-trivially on $V^{\ot r}$.
Actually, there is a minor complication: Dixmier's proof only applies if
$\dim V_{\Oz} \neq \dim V_{\Oo}$\,. But if $\dim V_{\Oz} = \dim V_{\Oo}$\,,
we can embed $V$ into $W = V \op \BC$\,, where $\BC$ is regarded as a trivial
$\fg$-module. Then his arguments apply to $W$, and the tensorial powers of
$W$ are isomorphic to direct sums of tensorial powers of $V$. This proves
the proposition. \vspace{1ex}

As at the beginning of this section, let $\cP$ be the dual of the embedding
of $U(\fg_{\Oz})$ in $U(\fg)$. We have
\begin{eqnarray*}
 \cP(t_{i\mu}) &\!=\!& \cP(t_{\mu i}) = 0 \,,\\[0.5ex]
 \cP(\Ot_{i\mu}) &\!=\!& \cP(\Ot_{\mu i}) = 0\,,
                             \quad 1 \leq i \leq m\,,\; m < \mu \leq m+n \,.
\end{eqnarray*}
Set 
 \[ \Ae=\cP(\A) \,. \]
Then $\Ae$ has a Peter-Weyl type basis in terms of the matrix elements of
irreducible finite-dimensional representations of
$\fsl(m) \oplus \fsl(n) \oplus \fgl(1)$. Thus it follows from the discussion
of the last section that there exists a unique normalized left integral
 \[ \int_0: \Ae \lra \BC \;, \]
which also turns out to be right invariant (see the Appendix). Denote by
$\nu(\Gamma) \in \Ac$ the image of $\Gamma$ under the natural embedding
$U(\fg) \lra \Ac$. Recalling Lemma \ref{lr}, we have

\begin{theorem} The linear form $\int = \nu(\Gamma) \cd \int_0\cP$ is a
non-trivial left and right integral \\ on $\A$\,.\end{theorem} 

To see that $\int$ is indeed non-trivial, we consider $\int\Theta\oTheta$\,,
where 
\begin{eqnarray*}
 \Theta_i &\!=\!& t_{i,\,m+n}t_{i,\,m+n-1} \cdots t_{i,\,m+1}\,,\\
 \oTheta_i &\!=\!& \Ot_{i,\,m+n}\Ot_{i,\,m+n-1} \cdots \Ot_{i,\,m+1}\,,\quad
                                       i = 1,\,2,\,\ldots,\,m \,,\\
 \Theta &\!=\!& \Theta_m\Theta_{m-1} \cdots \Theta_1\,,\\
 \oTheta &\!=\!& \oTheta_m\oTheta_{m-1} \cdots \oTheta_1 \,.
\end{eqnarray*}
We have
 \[ \int\!\Theta\oTheta
        \,=\, \left\langle\Theta\oTheta\,,\,\BE\,\BF\right\rangle\!
             \int_0\cP\left(\det(t_{\mu\nu})\det(\Ot_{\mu\nu})\right)^m. \]
As
 \[ \det(t_{\mu\nu})\det(\Ot_{\mu\nu})(u) = \ve(u)\,,\;\;\forall\;
                                                   u \in U(\fg_{\Oz})\,, \]
we immediately obtain
 \[ \int_0\cP\left(\det(t_{\mu\nu})\det(\Ot_{\mu\nu})\right)^m = 1 \,. \]
By induction we can show that
 \[ \langle\Theta\oTheta\,,\,\BE\,\BF\rangle = (-1)^{\frac{1}{2}mn(mn+1)}\,,\]
hence
 \[ \int\Theta\oTheta = (-1)^{\frac{1}{2}mn(mn+1)} \,. \]

\begin{example} The Lie supergroup $\mbox{\rm OSP}(2|2n)$ \end{example}
The Lie superalgebras $\fg=\fosp(2|2n)$ form the other series of type I
(basic classical) Lie superalgebras besides $\fsl(m|\,n)$. They share many
properties with the latter. In par\-ticular, the odd subspace of $\fosp(2|2n)$
is a direct sum of $\fg_{\Oo+}$ and $\fg_{\Oo-}$\,. Both $U(\fg_{\Oo+})$ and
$U(\fg_{\Oo-})$ are isomorphic to the Grassmann algebra on $2n$ generators.
The maximal even subalgebra of $\fosp(2|2n)$ is $\fsp(2n) \op \fgl(1)$, and
$\fg_{\Oz}$ and $\fg_{\Oo\pm}$ satisfy relations of the same form as
(\ref{same}).

The subspaces of $U(\fg_{\Oo\pm})$ of the highest Grassmann degree are both 
$1$-dimensional. We choose bases $\BE$ and $\BF$ for them, respectively, 
and set $\Gamma = \BE\,\BF$\,. Then 
 \[ X\Gamma \in U(\fg)\fg_{\Oz}\,,\;\; \forall\; X \in \fg \,. \]

Let $t$ be the defining representation of $\fosp(2|2n)$. It is known that
$t$ is self-dual. Introduce the matrix elements of $t$\,, 
 \[ t_{ab} \in U(\fg)^{\circ} \,, \quad a,b = 1,2,\ldots,2n+2 \,, \]
with $t_{ij}$ and $t_{\mu\nu}$ being even, and $t_{i\mu}$ and $t_{\mu i}$
odd, where $i,j = 1,2$\,; $\mu,\nu = 3,4,\ldots,2n+2$\,.

\begin{proposition} The elements $t_{ab}$ generate a sub-Hopf-superalgebra
$\A$ of $U(\fosp(2|2n))^{\circ}$, and $\A$ is dense in $U(\fosp(2|2n))^{\ast}$.
\end{proposition} 

{\em Proof}: This follows from the proof of Proposition \ref{Harish}.
							       \vspace{1ex}

Set $\Ae=\cP(\A)$. Then $\Ae$ admits a unique (up to scalar multiples) 
left integral $\int_0$\,. Denoting by $\nu(\Gamma)$ the canonical image of
$\Gamma$ in $\A^{\circ}$, we have 

\begin{theorem} The linear form $\int = \nu(\Gamma)\cd\int_0\cP$ is a
non-trivial left and right integral \\ on $\A$\,.
\end{theorem}

The $t_{i\mu}$ and $t_{\mu i}$ generate a Grassmann algebra contained in
$\A$\,. We take $\Theta$ to be a non-zero element of the highest degree 
in this Grassmann algebra. Then direct computations can show that 
 \[ \int\! \Theta \neq 0 \,.\vspace{2ex} \]

\noindent
{\bf Example 4}\hspace{1.0em} {\it The Lie supergroup $\mbox{\rm OSP}(1|2n)$} 
\nopagebreak \setcounter{example}{4}
\\[1.5ex]
Let us start with the simplest case, $n = 1$\,. The Dynkin diagram of
\nopagebreak
$\fosp(1|2)$ is just $\bullet$\,, and the simple Chevalley generators are
$\{e,\,f,\,h\}$\,, where $e$ and $f$ are odd while $h$ is even, with the
commutation relations
 \[ [h \,,\, e] = e \,,\quad [h \,,\, f] = -f \,,\quad [e \,,\, f] = h \,.\]
It is important to observe that $[e \,,\, e] = E$\,,\, $[f \,,\, f] = F$ and
$h$ span an $\fsl(2)$ subalgebra, which is the maximal even subalgebra
$\fosp(1|2)_{\Oz}$\,. This is a general feature of any type II superalgebra,
where some simple generators of the maximal even subalgebra are generated by
odd elements. We denote $\fg = \fosp(1|2)$,
$\fg_{\Oz} = \fsl(2) \subset \fosp(1|2)$, $\U = U(\fg)$ and
$\Ue = U(\fg_{\Oz})$. 

Now
 \[ \1 + ef + \U\fg_{\Oz} \]
is an invariant of the left $\U$-module $\U/\U\fg_{\Oz}$\,, and we have the
left and right integral
 \[ \int =\, \nu(\1+{ef}) \cd \int_0 \cP : \, \Uc \lra \BC \;, \]
where $\int_0 : \Uec \rightarrow \BC$ is the standard Haar functional on
$\Uec$\,. Consider $\int \1_{\Uc}$\,. We have 
 \[ \int \! \1_{\Uc}
       \,=\, \langle \1_{\Uc}, \1 + ef \rangle \! \int_0 \cP(\1_{\Uc})
                                        \,=\, \int_0 \1_{\Uec} \ne 0 \,. \]
That is, the integral does not vanish on the identity element of $\Uc$. It
follows from the discussion of Section 2 that all finite-dimensional
representations of $\fosp(1|2)$ are completely reducible, which, of course,
is a well-known fact.

The general case can be treated similarly. We do not go into details but
only mention that an even element $u_0 \in \U = U(\fosp(1|2n)$ such that
$\ve(u_0) \neq 0$ and such that $u_0 + \U \fg_{\Oz}$ is invariant in
$\U/ \U \fg_{\Oz}$ has been constructed by Djokovi\'c and Hochschild in
Ref.~\cite{DHo}. Moreover, they have proved the following theorem: \\[1.0ex]
Let $\fg$ be a finite-dimensional Lie superalgebra over a field of
characteristic zero. Then all finite-dimensional representations of $\fg$
are completely reducible if and only if the following two conditions are
satisfied. \\[1ex]
(1) The Lie algebra $\fg_{\Oz}$ is semisimple. \\[1ex]
(2) There is an element $u_0$ in $U(\fg)$ such that $u_0 + U(\fg)\fg_{\Oz}$
is an invariant element of $U(\fg)/ U(\fg)\fg_{\Oz}$ and satisfies
$\ve(u_0) \neq 0$\,. \\[1ex]
Visibly, in the cited reference the element $u_0$ has been a decisive tool
in the proof that all finite-dimensional representations of $\fosp(1|2n)$
are completely reducible. It is remarkable that in the present work it
serves to construct a left integral on $\Uc$ which does not vanish on the
unit element, a result which, in turn, implies the complete reducibility of
the $\Uc$-comodules and hence of the finite-dimensional $\U$-modules.

\begin{example} The Lie supergroup $\mbox{\rm OSP}(3|2)$ \end{example}
Let $\fg$ denote the Lie superalgebra $\fosp(3|2)$. It is the simplest of
those orthosymplectic Lie superalgebras which are not of type I and not one
of the special algebras $\fosp(1|2n)$. Its maximal even subalgebra is
$\fgz = \fso(3) \op \fsp(2)$. The $\fg$-module $U(\fg)/U(\fg)\fgz$ will be
denoted by $W$. We shall also need the quadratic Casimir element
$C \in U(\fg)$ and the corresponding Casimir operator $C_W$ acting on $W$.

In the subsequent investigation of the $\fg$-module $W$ we are going to use
the classification of finite-dimensional irreducible $\fg$-modules obtained
by Van der Jeugt in Ref. \cite{VdJ}. Both the $\fg$-modules and the
$\fgz$-modules are characterized by a pair of numbers
$p,q \in \{0,\oh,1,\th,\ldots\}$. By a slight abuse of notation, we
denote the corresponding
\vspace{0.3ex}
$\fg$-module by $[p,q]$, and the corresponding $\fgz$-module by $(p,q)$. (We
remark that $p$ is associated in the obvious way to $\fso(3)$ and $q$ to
$\fsp(2)$.)

A version of the Poincar\'e, Birkhoff, Witt theorem implies that $W$, regarded
as a $\fgz$-module, is isomorphic to the Grassmann algebra constructed over
$\fgo$\,. Using the representation theory of $\fsl(2)$, we conclude that the
$\fgz$-module $W$ decomposes into the direct sum of the modules contained in
the following list, where the first line gives the Grassmann degree to which
the modules underneath belong.
 \[ \begin{array}{c@{\hs}c@{\hs}c@{\hs}c@{\hs}c@{\hs}c@{\hs}c}
    0    &    1    &    2    &    3    &    4    &    5    &    6   \\[1.0ex]
  (0,0)  & (1,\oh) &  (1,1)  & (2,\oh) &  (1,1)  & (1,\oh) &  (0,0) \\[1.0ex]
         &         &  (2,0)  & (1,\oh) &  (2,0)  &         &        \\[1.0ex]
         &         &  (0,0)  & (0,\th) &  (0,0)  &         &
    \end{array} \]
Comparison with the lower-dimensional irreducible $\fg$-modules then shows
that for a Jordan-H\"older sequence of the $\fg$-module $W$ the irreducible
quotients must be isomorphic to the following modules:
 \[ [0,\th] \;,\;\; [1,1] \;,\;\; [1,\oh] \;,\;\; [0,0] \;,\;\; [0,0] \;. \]
For the convenience of the reader and for later use, we also note how these
modules decompose into irreducible $\fgz$-submodules, moreover, in the first
column we give the eigenvalue of the quadratic Casimir operator (normalized
as in Ref. \cite{VdJ}) in these modules.
 \[ \begin{array}{r@{\qquad}rcl}
  -6 & [0,\th] & \cong & (0,\th) \op (1,1) \op (1,\oh) \op (0,0) \\[1.0ex]
   0 &   [1,1] & \cong & (1,1) \op (1,\oh) \op (2,\oh) \op (2,0) \\[1.0ex]
   2 & [1,\oh] & \cong & (1,\oh) \op (2,0) \op (0,0)   \\[1.0ex]
   0 &   [0,0] & \cong & (0,0)
    \end{array} \]
Note that at this point it is obvious that the $\fg$-module $W$ is not
completely reducible: It is generated, as a $\fg$-module, by a
$\fgz$-invariant element, the multiplicity of $(0,0)$ in the $\fgz$-module
$W$ is equal to 4, but the length of the $\fg$-module $W$ (i.e., the number
of irreducible quotients of a Jordan-H\"older sequence) is equal to 5.

The eigenvalues given above imply that the primary decomposition of $W$ with
respect to $C_W$ takes the following form:
\begin{equation}
   W = W_{-6} \op W_2 \op W_0 \;, \label{decomp}
\end{equation}
where $W_r$\,, $r \in \{-6,2,0\}$, is the primary subspace of $W$ corresponding
to the eigenvalue $r$ of $C_W$. Of course, the $W_r$'s are $\fg$-submodules of
$W$. Regarded as $\fg$-modules, we have
 \[ W_{-6} \cong [0,\th] \quad,\quad W_2 \cong [1,\oh] \;, \]
whereas $W_0$ has a Jordan-H\"older sequence of the form
 \[ W_0 \supset W'_0 \supset W''_0 \supset \{0\} \;, \]
where one of the three modules $W_0/W'_0$\,, $W'_0/W''_0$\,, $W''_0$ is
isomorphic to $[1,1]$, while the other two are trivial one-dimensional. In
any case we have
\bea
 C_W(W_0) &\!\subset\!& W'_0 \;,   \label{cw1} \\[0.5ex]
 C_W(W'_0) &\!\subset\!& W''_0 \;, \label{cw2} \\[0.5ex]
 C_W(W''_0) &\!=\!& \{0\} \;.      \label{cw3}
\eea
We stress that while $W_{-6}$ and $W_2$ are eigenspaces of $C_W$, this is
not so for $W_0$\,. In fact, we shall see that the restriction of $C_W$ to
$W_0$ is not equal to zero but only nilpotent.

\begin{lemma}
The subspace $C_W(W_0)$ of $W_0$ is either a trivial one-dimensional
$\fg$-submodule of $W_0$ or else it is equal to $\{0\}$.
\end{lemma}

{\em Proof}: In the subsequent discussion, it is important to keep the
following fact in mind: \\[1ex]
$(\ast)$ The $\fg$-module $[1,1]$ does not contain a trivial $\fgz$-submodule.
                                                                      \\[1ex]
There are three cases to consider. \\[1ex]
a) The module $W_0/W'_0$ is isomorphic to $[1,1]$. \\[1ex]
This case is not possible since $W_0$\,, like $W$, is generated by a
$\fgz$-invariant element which, under the present assumption and because of
$(\ast)$, would have to belong to $W'_0$\,. \\[1ex]
b) The module $W'_0/W''_0$ is isomorphic to [1,1]. \\[1ex]
In this case, $W''_0$ consists of $\fg$-invariant elements, hence the
existence of non-zero $\fg$-invariant elements in $W_0$ is obvious.  However,
we want to find an explicit expression for them, and a first step towards
this end is the lemma. According to equation (\ref{cw2}) we have
$C_W(W'_0) \subset W''_0$\,. Using equation (\ref{cw3}) and recalling $(\ast)$,
we can even conclude that $C_W(W'_0) = \{0\}$. Thus $C_W$ induces a
$\fg$-module map $W_0/W'_0 \rightarrow W_0$\,, and this implies our claim.
Actually, it is easy to see that $C_W(W_0) \subset W''_0$\,. \\[1ex]
c) The module $W''_0$ is isomorphic to $[1,1]$. \\[1ex]
Equation (\ref{cw3}) says that $C_W(W''_0) = \{0\}$, hence $C_W$ induces a
$\fg$-module map $W'_0/W''_0 \rightarrow W'_0$ which, according to
equation (\ref{cw2}), is even a map into $W''_0$\,. Invoking $(\ast)$ we
conclude that $C_W(W'_0) = \{0\}$, and our claim follows as in part b).
\vspace{1ex} This proves the lemma.

Let us now recall the decomposition (\ref{decomp}) of $W$ and also the
fact that the $\fg$-module $W$ is generated by the element
$\1_{U(\fg)} + U(\fg)\fgz$\,. Then the lemma above can be rephrased as
follows: \\[1ex]
Either the element
 \[ z = C(C - 2)(C + 6) \in U(\fg) \]
belongs to $U(\fg)\fgz$\,, or else $z + U(\fg)\fgz$ is a non-trivial
invariant element of $U(\fg)/U(\fg)\fgz$\,.

Thus all that remains to be shown is that $z$ does not belong to
$U(\fg)\fgz$\,. This is an easy consequence of the Poincar\'e, Birkhoff,
Witt theorem, which allows us to construct a suitable basis of $W$. Actually,
the task can be simplified, as follows. The Casimir element $C$ can be
decomposed (in various ways) into the sum of two pieces,
 \[ C = C_o + C_e \;, \]
where $C_o$ is quadratic in the elements of $\fgo$\,, and where $C_e$
belongs to $U(\fgz)$. Since $C$ commutes with all elements of $U(\fg)$, it
follows that
 \[ z \in C_o(C_o - 2)(C_o + 6) + U(\fg)\fgz \;, \]
and hence we can replace $z$ by
 \[ z_o =  C_o(C_o - 2)(C_o + 6) \;. \]

Applying Theorem \ref{main} to $z$ or $z_o$ and recalling Lemma \ref{lr} we
obtain a non-zero left and right integral on $U(\fosp(3|2))^{\circ}$.

\section{Integrals on quantum supergroups}

We shall extend the construction of integrals on classical supergroups to
quantum supergroups at generic $q$\,. Recall that the Drinfeld-Jimbo quantum
superalgebra $U_q(\fg)$ associated with a simple basic classical Lie
superalgebra $\fg$ is usually defined with respect to the distinguished
simple root system of $\fg$ where only one odd simple root exists. By
removing the odd simple generators (but retaining the corresponding Cartan
generator), one obtains a graded quantum subalgebra
$U_q(\fg_0) \subset U_q(\fg)$, where $\fg_0 \subset \fg$ is an even
subalgebra of $\fg$, which is a reductive Lie algebra. We stress that while
for the basic classical  Lie superalgebras of type I we have
$\fg_0 = \fg_{\Oz}$\,, this is {\em not} the case for type II.

An important fact is that $U_q(\fg_0)$ forms a Hopf subalgebra of $U_q(\fg)$,
with its structure inherited from the latter. We have the following Hopf
superalgebra maps:
\[ \begin{array}{rl}
    \cI : \!\!\!& U_q(\fg_0) \lra U_q(\fg)\,,\\[0.8ex]
    \cP : \!\!\!& U_q(\fg)^{\circ} \lra U_q(\fg_0)^{\circ} \,,
\end{array} \]
where $\cI$ is the natural embedding and $\cP$ is induced from its dual
$\cI^{\ast}$.

A quantum supergroup associated with $\Uq$ is defined by specifying its
superalgebra of functions $\A$\,, where $\A$ should meet two basic
requirements, namely, it forms a sub-Hopf-superalgebra of $\Uq^{\circ}$, 
and it is dense in $\Uq^{\ast}$. In general, $\A$ is generated by the matrix
elements of some finite-dimensional irreducible representations of $\Uq$. The
structure of $\A$ associated with a type I quantum superalgebra has been
extensively studied. The fact that $\A$ is dense in $\Uq^{\ast}$ implies that
the natural Hopf superalgebra maps
\[ \begin{array}{rl}
     \nu : \!\!\!& \Uq \lra \Ac \,, \\[0.8ex]
     \hat{\cI} = \nu\cI : \!\!\!& \Uqo \lra \Ac \,,
    \end{array} \]
are embeddings. 

Denote $\Ae = \cP(\A)$. Then $\Ae$ separates points of $\Uqo$, i.e., it is
dense in $\Uqo^{\ast}$. Furthermore, $\Ae$ admits a Peter-Weyl type basis in
terms of the matrix elements of finite-dimensional irreducible representations
of $\Uqo$, and there exists a unique (up to scalar multiples) left integral 
 \[ \int_0:\, \A_e \lra \BC \,, \] 
which also turns out to be right invariant, and it is non-vanishing on
$\1_{\Ae}$\,.

Similar to the classical case, we consider
 \[ \int_0 \cP : \, \A \lra \BC \;, \]
which is clearly left invariant with respect to $\Uqo$, i.e.,
 \[ \hat{\cI}(u) \cd \int_0 \cP \,=\, \ve(u)\int_0 \cP\,,
                                           \quad\forall\; u \in \Uqo \,. \]
Let $K$ denote the ideal of $\Uqo$ defined by
 \[ K = \{u \in \Uqo |\, \ve(u) = 0 \} \,, \]
where $\ve$ is the counit of $\Uq$. Then
 \begin{equation} J = \Uq K \label{defJ} \end{equation}
is a left ideal of $\Uq$.

\begin{lemma}
If $\fg$ is one of the Lie superalgebras $\fsl(m|\,n)$ or $\fosp(2|2n)$
(i.e., if $\fg$ is basic classical of type I), the left ideal $J$ has finite
codimension in $\Uq$.
\end{lemma}

{\em Proof}: This follows immediately from the PBW theorems for these
quantum superalgebras established in \cite{slq, osp}. \vspace{1.5ex}

Clearly $\Uq/J$ forms a left $\Uq$-module under the natural action 
 \[ x(y + J) = xy + J \,, \quad \forall\; x,y \in \Uq \,.\]
Let $z + J$ be an invariant of $\Uq/J$, i.e.,
$x(z + J) = \ve(x)z + J \,,\; \forall\; x \in \Uq$. Non-trivial invariants of
this kind exist for type I quantum superalgebras, as we will see later.
However, we doubt that the type II quantum superalgebras admit such
invariants, as in this case $J$ is expected to have infinite codimension.

\begin{theorem} Let $\int = \nu(z) \cd \int_0 \cP$. Then $\int$ is a left
integral on $\A$\,, that does not depend on the representative of $z + J$.
As before, $\nu(z)$ is the image of $z$ under the natural embedding
$\Uq \lra \Ac$.
\end{theorem} 

{\em Proof}: The proof goes in the same way as in the classical case.

\begin{example} The quantum supergroup $\mbox{\rm SL}_q(m|\,n)$\end{example}
We study the quantum supergroup $\SLq$. The quantum superalgebra
$U_q(\fsl(m|\,n))$ is generated by the simple and the Cartan generators
 \[ E_{a,\,a+1}\,,\,\, E_{a+1,\,a}\,,\,\, k_a^{\pm 1}\,,\quad
                                         a = 1,\,2,\,\ldots,\,m+n-1 \,, \]
subject to the standard relations. (Here $k_a = K_a K_{a+1}^{-1}$ in the
notation of \cite{slq}.) The generators $E_{m,\,m+1}$ and $E_{m+1,\,m}$ are odd,
while all the others are even. Define recursively 
\begin{eqnarray*}
 E_{ab} &\!=\!& E_{ac}E_{cb} - q_c^{-1} E_{cb}E_{ac} \,,\\[0.5ex]
 E_{ba} &\!=\!& E_{bc}E_{ca} - q_c E_{ca}E_{bc} \,,\quad a < c < b \,,
\end{eqnarray*}
where $q_c = q^{(-1)^{[c]}}$. The vector representation $t$ of
$U_q(\fsl(m|\,n))$ is given by
 \[ t(E_{a,\,a \pm 1}) \,=\, e_{a,\,a \pm 1} \,, \]
 \[ t(k_a) \,=\, q_a^{e_{aa}}q_{a+1}^{-e_{a+1,\, a+1}}
        \,=\, 1 + (q_a - 1)e_{aa} + \left(q_{a+1}^{-1}
                             - 1\right)e_{a+1,\, a+1} \,. \vspace{0.7ex} \]
We shall denote the dual vector representation by $\Ot$, and let
 \[ t_{ab} \,,\, \Ot_{ab} \in U_q(\fsl(m|\,n))^{\circ}\,,\quad
                                          a,b = 1,\,2,\,\ldots,\,m+n \,, \]
be the matrix elements of $t$ and $\Ot$ respectively. Then the superalgebra
$\A$ of functions on the quantum supergroup $\SLq$ is defined to be the
subalgebra of $U_q(\fsl(m|\,n))^{\circ}$ generated by the $t_{ab}$\,,
$\Ot_{ab}$\,. In Ref. \cite{SLq} it was shown that

\begin{proposition}
The algebra $\A$ is a sub-Hopf-superalgebra of $U_q(\fsl(m|\,n))^{\circ}$ 
and is dense in $U_q(\fsl(m|\,n))^{\ast}$. 
\end{proposition}

The quantum even subalgebra $\Uqo$ is $U_q(\fsl(m) \op \fgl(1) \op \fsl(n))$
with generators
 \[ k_a^{\pm 1}\,,\;  E_{b,\,b+1}\,,\; E_{b+1,\,b}\,,
                         \quad a,b = 1,2,\ldots,m+n-1\,,\; b \neq m \,. \]
The images of $t$ and $\Ot$ under $\cP$ give rise to representations of
$\Uqo$, with 
 \[ \cP(t) = \left(\begin{array}{ll}
               \cP(t_{ij}) & 0 \\
                         0 & \cP(t_{\mu\nu})
             \end{array} \right)\;,  
\qquad 
  \cP(\Ot) = \left(\begin{array}{ll}
               \cP(\Ot_{ij})& 0 \\
                         0 & \cP(\Ot_{\mu\nu})
             \end{array} \right).  \]  
The matrix elements of these representations generate $\Ae$\,, which forms a
Hopf subalgebra of $\Uqo^{\circ}$. On $\Ae$ there exists a unique left integral
$\int_0$ which annihilates the matrix elements of all non-trivial irreducible
representations and satisfies
 \[ \int_0 \1_{\Ae} = 1 \,. \]

Introduce
\begin{eqnarray*}
 \BE_{\,i} &\!=\!& E_{i,\,m+1}\,E_{i,\,m+2}\,\cdots\,E_{i,\,m+n} \;,\\
 \BF_i &\!=\!& E_{m+n,\tty i}\,E_{m+n-1,\tty i}\,\cdots\,E_{m+1,\tty i} \;,\\
   \BE &\!=\!& \BE_{\,m}\,\BE_{\, m-1}\cdots\,\BE_{\,1} \,,\\
   \BF &\!=\!& \BF_1\,\BF_2\,\cdots\,\BF_m \,,\\
 \Gamma&=&\BE\,\BF \,.
\end{eqnarray*}

\begin{lemma} Let $J$ be defined as in (\ref{defJ}). Then the image
of\/ $\Gamma$ under the canonical map
$U_q(\fsl(m|\,n)) \lra U_q(\fsl(m|\,n))/J$ is an invariant.
\end{lemma}

{\em  Proof}: In \cite{slq} it was shown that
\begin{eqnarray*}
 k_a \Gamma &\!=\!& \Gamma\, k_a \,, \quad\forall\; a \,,\\
 {} [E_{c,\,c+1},\,\BE\,] &\!=\!& [E_{c,\,c+1}\,,\,\BF\,] = 0\,,\quad
                                                            c \neq m \,,\\
 {} [E_{c+1,\,c}\,,\,\BE\,] &\!=\!& [E_{c+1,\,c}\,,\,\BF\,] = 0\,,\quad
                                                            c \neq m \,.
\end{eqnarray*}
It is also clear that
 \[ E_{m,\,m+1}\Gamma = 0 \,. \]
This immediately leads to
 \[ E_{i,\,m+1} \Gamma = 0 \,, \quad\forall\; i \leq m \,. \]
What remains to be shown is that
\begin{equation}
  E_{m+1,\,m}\Gamma \in J \,.\label{master}
  \end{equation}
By using the fact that $E_{m+1,\,m}$ $q$-anticommutes with all
$E_{\mu, i}$\,, $\mu \geq m+1$, $i\leq m$\,, and $(E_{m+1,\,m})^2 = 0$\,,
we have
 \[ E_{m+1,\,m}\,\BF = 0 \,. \]
Thus
 \[ E_{m+1,\,m}\Gamma = [E_{m+1,\,m}\,,\,\BE\,]\,\BF\,. \]
To determine the right hand side, we need the following commutation relations
\begin{eqnarray*}
 [E_{m+1,\,m}\,,\,\BE_{\,i}]
      &\!=\!& q^{m+n-2}E_{i,\,m+2}\,E_{i,\,m+3}\,
                  \cdots\,E_{i,\,m+n}k_m E_{i,\,m} \,,\; i < m \,,\\[0.5ex]
 [E_{i,\,m}\,,\,\BE_{\,j}] &\!=\!& 0 \,,\quad i > j \,.
\end{eqnarray*}
Now
 \[ [E_{m+1,\,m}\,,\,\Gamma] = [E_{m+1,\,m}\,,\,\BE_{\,m}]
                              \,\BE_{\, m-1}\,\cdots\,\BE_{\,1}\,\BF \;, \]
where $[E_{m+1,\,m}\,,\,\BE_{\,m}]$ can be easily calculated to yield
\vspace{-1.0ex}
 \[ \begin{array}{l}
   [E_{m+1,\,m}\,,\,\BE_{\,m}] \,=\, \ds{\frac{k_m-k_m^{-1}}{q-q^{-1}}\,
  E_{m,\,m+2}\,\cdots\, E_{m,\,m+n} } \\[2.0ex]
   \ds{\hspace{7.5em} + \sum_{\alpha=2}^n (-1)^\alpha q^{-(n-\alpha)}
           E_{m,\,m+1}\,\cdots\,\widehat{E}_{m,\,m+\alpha}
                \,\cdots\,E_{m,\,m+n}E_{m+1,\,m+\alpha} \,k_\alpha^{-1}} \,,
                                                            \vspace{-1.0ex}
\end{array} \]
with $\widehat{E}_{m,\,m+\alpha}$ indicating that $E_{m,\,m+\alpha}$ is
removed from the second term. By using
 \[ E_{m+1,\,m+\alpha}\,\BE_{\,i}
           - q^{-2}\,\BE_{\,i}\,E_{m+1,\,m+\alpha} = 0\,,
                      \quad i = 1,2,\ldots,m\,,\; \alpha = 2,3,\ldots,n \,,\]
we immediately see that (\ref{master}) indeed holds. \vspace{0.5ex}

Let $\nu : \Uq \rightarrow \A^{\circ}$ be the natural embedding.

\begin{theorem} There exists the following non-trivial left integral on
${\rm SL}_q(m|\,n)$\,:
 \[ \int = \nu({\Gamma}) \cd \int_0 \cP \,. \]
\end{theorem}

\begin{example} The quantum supergroup ${\rm OSP}_q(2|2n)$ \end{example}
We denote  by $\fg$ the Lie superalgebra $\fosp(2|2n)$ and recall that in this
case $\fg_0 = \fg_{\Oz}$ is the maximal even subalgebra $\fsp(2n)\oplus\fgl(1)$
of $\fg$\,. Introduce the $(n+1)$-dimensional Minkow\-ski space $\fh^{\ast}$
with a basis $\{\delta_i \,|\, i=0, 1, 2,\ldots,n\}$ and the bilinear form 
$(\ ,\ ): \fh^{\ast}\times \fh^{\ast}$ $\rightarrow \BC$ defined by 
 \[ (\delta_i\tty,\tty\delta_j)
  = -(-1)^{\delta_{0,i}} \, \delta_{i,\tty j} \;, \quad \forall\; i,j \,. \]
Then the simple roots can be expressed as 
$\alpha_i = \delta_i-\delta_{i+1}$\,, $0 \leq i < n$\,, 
$\alpha_n = 2\delta_n$\,, with $\alpha_0$ being the unique odd simple root.
A convenient version of the  Cartan matrix $A = (a_{ij})_{i,j=0}^n$ is  
$a_{ij} = 2(\alpha_i, \alpha_j)/(\alpha_i, \alpha_i)$, $\forall\: i>0$\,,  
$a_{0,j}=(\alpha_0, \alpha_j)$. The quantum superalgebra $\Uq$ is the
universal complex superalgebra with generators
$\{ k_i^{\pm 1},\, e_i\,,\, f_i\,, \; i\in\Nn\}$, $\Nn = \{0,1,2,\ldots,n\}$,
where $e_0$ and $f_0$ are odd and the rest are even. The defining relations
are
\bea 
& &k_i k_j = k_j k_i \;,
\quad k_i k_i^{-1} = k_i^{-1} k_i = \1 \,,\nonumber \\[0.5ex] 
& &k_i e_j k_i^{-1}= q_i^{a_{ij}/2}e_{j} \,, \quad 
   k_i f_j k_i^{-1}=q_i^{-a_{ij}/2}f_{j} \,, \nonumber\\[1.0ex] 
& &{[}e_i, f_j] = \delta_{ij}(k_i^2-k_i^{-2})/(q_i-q_i^{-1})\,,
\quad i,j \in \Nn \,, \nonumber\\[1.0ex]
& & (e_0)^2=(f_0)^2=0 \,,\nonumber \\[0.5ex] 
&&\sum_{\mu=0}^{1-a_{ij}}(-1)^\mu 
\left[\begin{array}{c}      
       1-a_{ij}\\  
      \mu 
      \end{array} 
\right]_{q_{i}} 
e_i^{1-a_{ij}-\mu }e_je_i^\mu =0 \,, \;\; i \neq 0 \,, \nonumber \\   
& &\sum_{\mu=0}^{1-a_{ij}}(-1)^\mu 
\left[\begin{array}{c}      
       1-a_{ij}\\  
      \mu 
      \end{array} 
\right]_{q_{i}} 
f_i^{1-a_{ij}-\mu }f_jf_i^\mu =0 \,, \;\;  i\ne 0 \,, \nonumber             
\eea 
where $\Big[\begin{array}{c} \!m\! \\ n \end{array}\Big]_q$ 
is a $q$-binomial coefficient. As is well-known, \vspace{-0.3ex}
the quantum superalgebra
$\Uq$ has the structure of a Hopf superalgebra. Note that \vspace{0.4ex}
$\{e_i \,,\, f_i \,,\, k_i^{\pm 1} \,|\, i = 1,2,\ldots,n\}$ generate a Hopf
subalgebra $U_q(\fsp(2n)) \subset \Uq$. Together with $\{k_0^{\pm 1}\}$,
\vspace{0.1ex}
they  generate $\Uqo = U_q(\fsp(2n)\oplus \fgl(1))$. \vspace{0.2ex}

Define  the odd elements  
\bea
\psi_1&=&e_0 \,,  \nonumber \\ 
\psi_{i+1}&=&\psi_{i}e_i-q e_i\psi_{i} \,, \ \  1\leq i<n \,, \nonumber \\ 
\psi_{- n}&=&\psi_n e_n - q^2 e_n \psi_n \,,  \nonumber \\ 
\psi_{- i}&=&\psi_{-i-1}e_i-q e_i\psi_{-i-1} \,, \ \ 1\leq i<n \,;
                                                              \nonumber \\
\phi_0&=&f_0 \,,  \nonumber \\
\phi_{i+1}&=& f_i \phi_i - q^{-1} \phi_i f_i \,, \ \ 1\leq i<n \,,\nonumber \\
\phi_{- n}&=&f_n \phi_n - q^{-2} \phi_n f_n \,, \nonumber \\
\phi_{- i}&=& f_i \phi_{-i-1} -q^{-1} \phi_{-i-1} f_i \,, \ \ 1\leq i<n \,,
						                    \nonumber
\eea 
which satisfy the following relations
\bea 
\psi_{\pm i} \psi_{\pm j} +q^{\pm 1}\psi_{\pm j} \psi_{\pm i}&=&0 \,, 
\ \  i\leq j \,, \nonumber \\   
\psi_i\psi_{- j} + q \psi_{- j}\psi_{i}&=&0 \,, \ \ \forall\; i\neq j \,, 
                                                              \nonumber \\  
\psi_n\psi_{- n}+q^2\psi_{- n}\psi_n&=&0 \,,   \nonumber \\ 
\psi_{-i-1}\psi_{i+1}+\psi_{i+1} \psi_{-i-1} 
+q\psi_{- i}\psi_i+q^{-1}\psi_i\psi_{- i} &=&0 \,, \ \ i<n \,; \nonumber \\
\psi_j e_i - q^{(\alpha_i, \delta_0 -\delta_j)}e_i\psi_j
&=&\delta_{i j}\psi_{i+1} \,, \ \ \forall\; i,j \,, \nonumber \\  
\psi_{-j} e_i - q^{(\alpha_i, \delta_0 +\delta_j)}e_i\psi_{-j}
&=&\delta_{i+1, j}\psi_{-i+1} \,, \ \ i>1 \,, \nonumber
\eea 
and also similar relations for $\phi_{\pm i}$\,, where $\psi_{n+1}$ and 
$\phi_{n+1}$ are understood as $\psi_{-n}$ and $\phi_{-n}$ respectively.
Let \vspace{-1.0ex} 
\be
E_{1,\,2}&=&e_1 \,,\\
E_{1,\,i+1}&=&E_{1,\,i} e_i - q e_i E_{1,\,i} \,, \quad 1<i<n \,,\\
E_{1,\,\ol{n}}&=& E_{1,\,n} e_n -q^2 e_n E_{1,\,n} \,,\\
E_{1,\,\ol{i}}&=& E_{1,\,\ol{i+1}} e_i - q e_i E_{1,\,\ol{i+1}} \,, 
                                                    \quad 1<i<n \,, \\
E_{1,\,\ol{1}}&=& E_{1,\,\ol{2}} e_1 q^{-1} - q e_1 E_{1,\,\ol{2}} \,,
\ee
where we have introduced the notation $\ol{i} = -i$\,. Then 
\bea
\{ \psi_i \,,\, f_0\}&\!=\!& E_{1,\,i}\, k_0^{-2} \,,\nonumber \\
\{\psi_{-i} \,,\, f_0\}&\!=\!& E_{1,\,\ol{i}}\, k_0^{-2} \,. \nonumber
\vspace{-1.0ex}
\eea
Define \vspace{-1.5ex}
\bea 
 \BE &\!=\!& \psi_1 \psi_2 \cdots \psi_n
                     \psi_{-n} \psi_{-n+1} \cdots \psi_{-1} \,, \nonumber \\ 
 \BF &\!=\!& \phi_{-1}\phi_{-2} \cdots \phi_{-n} 
                     \phi_n\phi_{n-1} \cdots \phi_1 \,, \nonumber\\ 
 \Gamma &\!=\!& \BE\,\BF \,. \nonumber   
\vspace{-0.5ex}
\eea
We have \vspace{1ex} 

\begin{lemma}
Let $J$ be defined as in (\ref{defJ}). Then
 \[ \begin{array}{rl}
  i)& [v\,,\,\BE\,] = [v\,,\, \BF\,] = 0 \,,\quad \forall\;
                                v \in U_q(\fsp(2n))\subset\Uqo \,,\\[0.7ex]
 ii)& [u \,,\,\Gamma] = 0 \,,\quad \forall\; u\in \Uqo \,, \\[0.7ex]
iii)& x\Gamma \in \ve(x)\Gamma + J \,, \quad \forall\; x\in\Uq \,.  
    \end{array} \]
\end{lemma}

Of particular importance for us is the vector representation $t$ of $\Uq$.
Introduce the index $a=i$ or $\ol{i}$, with $i = 0,1,\ldots,n$\,,
$\,\ol{i} =\ol{0},\ol{1},\ldots,\ol{n}$\,. We have \\[-3.0ex]
\begin{minipage}[t]{6.25in}
 \[ t(e_0) = e_{0,1} + e_{\ol{1},\,\ol{0}} \;,\quad
               t(f_0) = e_{1,\,0} - e_{\ol{0},\ol{1}} \;, \vspace{-0.8ex}\]
 \[ t(e_i) = e_{i,\,i+1} - e_{\ol{i+1},\,\ol{i}} \;,\quad
  t(f_i) = e_{i+1,\,i} - e_{\ol{i},\,\ol{i+1}} \;,\quad 1 \leq i < n \,, \]
 \[ t(e_n) = e_{n,\,\ol{n}} \;, \quad t(f_n) = e_{\ol{n},\,n} \;,
                                                          \vspace{0.3ex} \]
 \[ t(k_i) = q_i^{H_i/2} \,,\quad 0 \leq i \leq n \,, \]
\end{minipage}
where
\be
 H_0 & \!=\!& \delta^{\ast}_0 + \delta^{\ast}_1 \;, \\
 H_i & \!=\!& \delta^{\ast}_i - \delta^{\ast}_{i+1} \;,\quad  0<i<n \,, \\
 H_n & \!=\!& \delta^{\ast}_n \,; \\
 \delta^{\ast}_i &\!=\!& e_{i,\,i} - e_{\ol{i},\,\ol{i}} \;, \quad
                                                       0 \leq i \leq n \,.
\ee

Let $t_{ab} \in \Uq^{\circ}$,
$a,b = 0,1,\ldots,n,\ol{0},\ol{1},\ldots,\ol{n}$ be the matrix elements of
the vector representation $t$\,,   
 \[ \langle t_{ab} \,,\, x \rangle \,=\, t(x)_{ab} \,, \quad
                                                \forall\; x \in \Uq \,. \]
We will take the algebra $\A$ of functions on $\OSPq$ to be the subalgebra
of $\Uq^{\circ}$ generated by the elements $t_{ab}$\,. In \cite{OSP} we have
shown: 

\begin{proposition}
The algebra $\A$ is a sub-Hopf-superalgebra of $\U_q(\fosp(2|2n))^{\circ}$
and is dense in $\U_q(\fosp(2|2n))^{\ast}$.
\end{proposition}

As usual, let $\cP : \Uq^{\circ} \rightarrow \Uqo^{\circ}$ be the map
induced by the dual of the embedding $\cI : \Uqo\rightarrow\Uq$, let
$\nu : \Uq \rightarrow \A^{\circ}$ be the canonical map, and let
$\hat{\cI} = \nu\,\cI$. Set $\Ae=\cP(\A)$. Then $\Ae$ admits a left integral
$\int_0$\,, which we normalize by setting $\int_0 \1_{\Ae} = 1$\,. Now 
 \[ \int_0\cP : \A \lra \BC \]
is a well-defined linear map, which is left invariant with respect to
$\Uqo\subset\Uq\;$:
 \[ \hat\cI(u) \cd \int_0 \cP \,=\, \ve(u)\int_0 \cP \,,
                                            \quad\forall\; u \in \Uqo \,.\]
We define 
\bea
 \int = \nu({\Gamma})\cdot \int_0\cP \,. \nonumber
\eea

\begin{theorem}
The linear form $\int\!:\tty \A \rightarrow \BC$ is a left integral on $\A$\,.
\end{theorem}

Consider $\int\Lambda$\,, where 
 \[ \Lambda \,=\, t_{\ol{1}\,\ol{0}} \cdots t_{\ol{n}\,\ol{0}}\,
                       t_{n\,\ol{0}} \cdots t_{1\,\ol{0}}\,
                            t_{\ol{1}\,0} \cdots t_{\ol{n}\,0}\,
	                         t_{n\,0} \cdots t_{1\,0} \;. \]
Using the following property of the Hopf superalgebra homomorphism $\cP\;$:
\be 
 \cP(t_{a\,0}) = \cP(t_{a\,\ol{0}}) = 0 \,, && \forall\; a \neq 0,\ol{0} \,,\\
 \cP(t_{0\,\ol{0}}) = \cP(t_{\ol{0}\,0}) = 0 \,, 
\ee 
we have 
 \[ \int\!\Lambda \,=\, \langle \Lambda\,,\, \Gamma \rangle 
    \int_0 \cP\left((t_{\ol{0}\,\ol{0}})^{2n} (t_{0\,0})^{2n}\right) \,. \]
Now 
 \[ \cP(t_{\ol{0}\,\ol{0}}\,t_{0\,0}) = \cP(t_{0\,0}\,t_{\ol{0}\,\ol{0}})
                                           = \1_{U_q(\fg_0)^{\circ}} \,, \]
thus
 \[ \int\!\Lambda = \langle \Lambda\,,\, \Gamma \rangle \,, \] 
which does not vanish if its $q \rightarrow 1$ limit is non-zero. A brute
force calculation shows 
 \[ |\langle \Lambda \,,\, \Gamma \rangle| \rightarrow 1,
                     \;\;\mbox{as}\;\; q \rightarrow 1 \,. \vspace{3.0ex}\]

\section{Discussion}

In the present work we have introduced and investigated the integrals on
Hopf superalgebras, with special emphasis on the classical and quantum
supergroups. In the undeformed case, there is obviously one problem that we
have not been able to solve completely, namely, to prove the existence of
non-zero integrals for all of the basic classical Lie supergroups. The
reason is that we have not really understood the structure of the induced
$\U(\fg)$-module $W = \U(\fg)/\U(\fg)\tty\fgz$\,. It is tempting to conjecture
that the invariant in $W$ we are looking for can always be expressed in terms
of Casimir elements (or even in terms of the quadratic Casimir element) as in
the example of ${\rm OSP}(3|2)$.

In the quantum case, we have only been able to treat the type I supergroups.
In particular, we could not say anything about the orthosymplectic quantum
supergroups. There are clear indications that our method will not work (or,
at least, has to be modified) in this case. However, one should remember
that, at present, only very little is known about the orthosymplectic quantum
supergroups anyway.
\vspace{3ex}

\noindent
{\bf Acknowledgement} \\
The present work was started during a visit of the first--named author to
the Department of Pure Mathematics at the University of Adelaide; it was
continued when both authors attended the seminar and workshop on
Cooperative Phenomena in Statistical Physics, held at the Max Planck
Institute for Physics of Complex Systems in Dresden; and it was finished
during a visit of the first--named author to the Department of Mathematics
at the University of Queensland, Brisbane. The kind invitations to these
institutes and the hospitality extended there are gratefully acknowledged.
\vspace{4.0ex}

\noindent
{\LARGE\bf Appendix\\[-3ex]}

\begin{appendix}

\section{Description of $U(\fgl(1))^{\circ}$}
In Example 2 of Section 3 we need to choose a (left) integral on
 \[ \Ae \subset \Uec \;, \vspace{-2ex} \]
where
 \[ \Ue = U(\fgz) = U(\fsl(m) \op \fsl(n) \op \fg(1))
                      \cong U(\fsl(m)) \ot U(\fsl(n)) \ot U(\fgl(1)) \;, \]
and hence
 \[ \Uec \cong U(\fsl(m))^{\circ} \ot U(\fsl(n))^{\circ}
                                                  \ot U(\fgl(1))^{\circ} \]
(the isomorphisms are to be interpreted in the Hopf algebra sense). According
to the discussion in Section 2, the Hopf algebras $U(\fsl(n))^{\circ}$ are
sufficiently well understood. In particular, there is a unique (up to scalar
multiples) left integral on $U(\fsl(n))^{\circ}$, which turns out to be right
invariant as well. For $\fosp(2|2n)$ and for the quantum counter parts the
situation is similar. Correspondingly, in the present appendix we would like
to comment on $U(\fgl(1))^{\circ}$. Needless to say, the results to be
presented are well-known \cite{PTa, CPi}, and we summarize them here in order
to clarify some slightly subtle issues.

The Lie algebra $\fgl(1)$ is one-dimensional, hence $U(\fgl(1))$ is
isomorphic (as a Hopf algebra) to the polynomial algebra $\BC[X]$ in one
indeterminate $X$. The Hopf algebra structure is the one known from
enveloping algebras: The structure maps are uniquely fixed by the equations
 \[ \begin{array}{rcl}
      \Delta(X) \!&=&\! X \ot \1 + \1 \ot X     \\[0.5ex]
         \ve(X) \!&=&\! 0                     \\[0.5ex]
           S(X) \!&=&\! - X \;.
    \end{array} \]
It follows that
 \[ \Delta(X^r) \,=\, \sum_{s=0}^r \Big(\begin{array}{c}
					  \! r \! \\ \! s \!
                                       \end{array} \Big) \, X^s \ot X^{r-s}
                                                                      \;, \]
for all integers $r \geq 0$\,.

The finite dual $\BC[X]^{\circ}$ of $\BC[X]$ can be described as follows.
Define, for any element $a \in \BC$ and any integer $r \geq 0$\,, the
linear form $u^r_a$ on $\BC[X]$ by
 \[ \langle u^r_a\,, P \rangle
 = \frac{d^{\, r} \! P}{d X^{r\st{-0.7ex}{0.1ex}}} \Big|_{X=a}
                                     \;, \quad\forall\; P \in \BC[X] \;. \]
Using some elementary algebra, it is not difficult to prove that these linear
forms, with $a$ and $r$ as described above, form a basis of the vector space
$\BC[X]^{\circ}$. The multiplication in $\BC[X]^{\circ}$ is given by
 \[ u^r_a u^s_b = u^{r+s}_{a+b} \;, \]
in particular, the unit element is equal to $u^0_0$ (which is the counit of
$\BC[X]$), the coproduct is given by
 \[ \Delta(u^r_a) \,=\, \sum_{s=0}^r
                        \Big(\begin{array}{c}
                         \! r \! \\ \! s \!
                        \end{array} \Big) \, u^s_a \ot u^{r-s}_a \;, \]
the counit by
 \[ \ve(u^r_a) = \delta_{r,\, 0} \;, \]
and the antipode by
 \[ S(u^r_a) = (-1)^r u^r_{-a} \;, \]
where, in all cases, $a,b \in \BC$ and $r,s \geq 0$ are integers.

Let us next recall that the dual $\BC[X]^{\ast}$ of the vector space $\BC[X]$
can be identified (in various ways) with the space of formal power
series $\BC[[Y]]$ in one indeterminate $Y$. If the dual pairing
 \[ \langle \,\:,\; \rangle : \BC[[Y]] \times \BC[X] \lra \BC \]
is chosen such that
 \[ \Big\langle \sum_{n \geq 0} c_n Y^n, X^r \Big\rangle
                                \,=\, r! \, c_r \;, \quad\forall\; r \;, \] 
then the coalgebra structure of $\BC[X]$ induces just the usual algebra
structure on $\BC[[Y]]$. Using this identification, the corresponding
injection
 \[ \BC[X]^{\circ} \lra \BC[[Y]] \]
is given by
 \[ u^r_a \lra Y^r \!\!\:\exp(aY) \;, \]
which immediately gives the product rule for the $u^r_a$\,'s. Similarly, we
find
 \[ \Delta(Y^r \!\!\:\exp(aY)) \,=\, (Y \ot \1 + \1 \ot Y)^r
                                           (\exp(a Y) \ot \exp(a Y)) \;. \]
Under the canonical embedding of $\BC[[Y]] \ot \BC[[Y]]$ into
$\BC[[Y \ot \1,\1 \ot Y]]$\,, the algebra of formal power series in $Y \ot \1$
and $\1 \ot Y$, the right hand side of this equation can be written in the form
 \[ (Y \ot \1 + \1 \ot Y)^r \exp(a(Y \ot \1 + \1 \ot Y)) \;. \] 
In this sense, the coproduct in $\BC[X]^{\circ}$ is fixed by the simple rule
 \[ \Delta (Y) = Y \ot \1 + \1 \ot Y \,, \]
just as for $\BC[X]$.

Let us now turn to the object of our main concern, the integrals. It is easy
to see that on $\BC[X]^{\circ}$ a non-trivial integral does not exist.
However, there is a way out. Obviously, the elements $u^0_a$\,, $a \in \BC$\,,
span a Hopf subalgebra $\K$ of $\BC[X]^{\circ}$, and the linear form
$\int$ on $\K$\,, defined by
 \[ \int\! u^0_a = \delta_{a,\,0} \;, \quad\forall\; a \in \BC \;, \]
is a left and right integral on $\K$\,. Note that the $u^0_a$\,'s are exactly
the characters of the algebra $\BC[X]$, i.e., the group-like elements of 
$\BC[X]^{\circ}$, and that $\K$ is isomorphic to the group Hopf algebra of
the additive group $\BC$\,.

Now we recall that, for an arbitrary algebra $A$ (associative, with unit
element), the finite dual $A^{\circ}$ consists exactly of the matrix elements
(regarded as linear forms on $A$) of the representations of $A$\,. (Here and
in the following, all representations are assumed to be finite-dimensional.)
It is easy to see that the matrix elements of the completely reducible
representations of $\BC[X]$ (i.e., the representations for which the image
of $X$ is diagonalizable) belong to $\K$\,, whereas the other elements of
$\BC[X]^{\circ}$ stem from those representations which are not completely
reducible.  Note that, once again, the close relationship between complete
reducibility and the existence of non-trivial integrals shows up.

Returning to the situation at the beginning of this appendix, we have to
assume that
 \[ \Ae \subset U(\fsl(m))^{\circ} \ot U(\fsl(n))^{\circ} \ot \K \;. \]
According to the foregoing discussion, this corresponds to the requirement
to consider only those representations of $\fg = \fsl(m|\, n)$ for which the
one-dimensional center of $\fgz$ is represented by diagonalizable operators,
which is usually assumed anyway.

\end{appendix}

\end{document}